\documentstyle[12pt]{article}
\begin{document}
\title{Morita equivalences between fixed point algebras and crossed 
products}
\author{Chi-Keung Ng\footnotemark[1]}
\addtocounter{footnote}{1}\footnotetext{This work is supported by 
the Croucher Foundation.}
\date{}
\maketitle
\begin{abstract}
In this paper, we will prove that if $A$ is a $C^*$-algebra with an 
effective coaction $\epsilon$ by a 
compact quantum group, then 
the fixed point algebra and the reduced 
crossed product are Morita equivalent. 
As an application, we prove an imprimitivity type theorem for crossed 
products of coactions by discrete Kac $C^*$-algebras. 
\end{abstract}

\par\medskip\par\medskip
\noindent {\bf 0. Introduction}

\par\medskip\par\medskip
After proving an imprimitivity type theorem for multiplicative 
unitaries of discrete type (actually, those come from discrete Kac algebras) 
in [7, 5.2], we noticed that the same method can be use to prove the fact 
that the fixed point algebra and the reduced crossed product of an effective 
coaction (see Definition 2.2) on a $C^*$-algebra A by a Woronowicz 
$C^*$-algebra (i.e. a compact quantum group)
are Morita equivalent (see Theorem 2.13). 
However, since $A$ may not be unital, we need a slightly more 
general version of Watatani's $C^*$-basic construction. We do 
this in Section 1. 

\par\medskip
In Section 2, we prove the main theorem. We first prove 
this in the case when the coaction is injective and the Hopf $C^*$-algebra 
is the reduced one. We then show how we can relax the condition to the case 
when the coaction is not injective and the Hopf $C^*$-algebra may be 
the full one. 

\par\medskip
In Section 3, we give some applications of the main theorem. 
In particular, we show that we don't need to assume amenability in the 
imprimitivity theorem in [7, 5.2] (see [7, 5.7(a)]). 
Moreover, we also prove an imprimitivity type theorem 
for crossed products of coactions by discrete type Kac $C^*$-algebras. 
Finally, we give an application to the case of discrete group coactions.

\par\medskip\par\medskip

\noindent {\bf 1. $C^*$-basic constructions for non-unital $C^*$-algebras}

\par\medskip\par\medskip
We would like to start with a slightly more general version of Watatani's 
$C^*$-basic construction (see [14]) in the case when the $C^*$-algebras 
may not be unital. 
Actually, it is just an application of the original version by Watatani. 
Let $B$ be a $C^*$-algebra and let $A$ be a 
$C^*$-subalgebra of $B$ that contains an approximate unit of $B$. 
Let $E$ be a faithful conditional expectation from $B$ to $A$. 
In this case, we have the following easy lemma.

\par\medskip
\noindent {\bf Lemma 1.1:} $A$ has a unit if and only if $B$ has. Moreover, if 
it is the case, then the unit of $B$ is in $A$.

\par\medskip
As in [14, 2.1], we define an A-valued inner product on $B$ 
by $\langle b,c\rangle  = E(b^*c)$. 
Since $E$ is a faithful conditional expectation, the inner product 
$\langle \cdot,\cdot\rangle $ makes $B$ into a pre-Hilbert $A$-module. 
Let $\cal F$ be the completion of $B$ with respect to this inner product 
and let $\eta$ be the canonical map from $B$ to $\cal F$. 
Consider a *-homomorphism $\lambda$ from $B$ to ${\cal L(F)}$ defined by 
$\lambda(a)(\eta(b)) = \eta(ab)$ and 
an element $e_A\in{\cal L(F)}$ given by $e_A(\eta(b))=\eta(E(b))$. 
Using the same argument as in [14, 2.1], we can show that they are 
well defined and $e_A$ is clearly a projection in $\cal L(F)$. 
As in [14, 2.1.2], we define the $C^*$-basic construction as follows:

\par\medskip
\noindent {\bf Definition 1.2:} Let $C^*\langle B,e_A\rangle $ be the closed 
linear span of the set $\{ \lambda (a) e_A\lambda(b): a, b\in B \}$ 
and call it the (reduced) $C^*$-basic construction. 

\par\medskip
Note that in [14], there are two kinds of $C^*$-basic construction but since 
they are isomorphic, we only consider one of them here. 
It is clear that $C^*\langle B,e_A\rangle  = {\cal K(F)}$ (where 
$\cal K(F)$ is the closed linear span of the finite rank operators 
on $\cal F$). 
In the case when the $C^*$-algebras are unital, we have the following result 
from [14, 2.2.11]:

\par\medskip
\noindent {\bf Theorem 1.3:} (Watatani) Let B be a unital 
$C^*$-algebra and let $A$ be a $C^*$-subalgebra of $B$ that contains 
the unit of $B$. 
Let $E$ be a faithful conditional expectation 
from $B$ to $A$. 
If $B$ acts on a Hilbert space $H$ faithfully (and non-degenerately) 
and $e$ is a projection on $H$ such that 
\par
\noindent (i) $ebe = E(b)e$ for all $b\in B$ and 
\par
\noindent (ii) the map that sends $a\in A$ to $ae\in {\cal L}(H)$ 
is injective, 
\par
\noindent then the norm closure of $BeB$ is isomorphic to 
$C^*\langle B,e_A\rangle $ canonically.

\par\medskip
We would like to give a similar result in the non-unital case. 
In the remainder of this section, we assume that $B$ is non-unital. 
Let $B^1=B\oplus \bf C\cdot 1$ be the unitalisation of $B$ and let 
$A^1$ be the $C^*$-subalgebra of $B^1$ generated by $A$ and $1$.
Define a map $E^1$ from $B^1$ to $A^1$ by $E^1(x+\mu 1) = E(x)+\mu 1$. 
Then we have:

\par\medskip
\noindent {\bf Proposition 1.4:} 
$E^1$ is a faithful conditional expectation from $B^1$ to $A^1$. 
\par
\noindent {\bf Proof:} We first show that $\parallel E^1\parallel = 1$. 
Note that $A^1 \subseteq M(A) \subseteq M(B)$ (since $A$ contains an 
approximate unit of $B$). Then $\parallel E^1(x+\mu1)\parallel_{M(A)} = 
sup \{\parallel E^1(x+\mu 1)y\parallel: y\in A, \parallel y\parallel 
\leq 1 \}=$ $sup \{\parallel E(xy+\mu y)\parallel: y\in A, \parallel 
y\parallel \leq 1 \} \leq$ $sup \{\parallel xy+\mu y\parallel: y\in A, 
\parallel y\parallel \leq 1 \} \leq$ $\parallel x+\mu1\parallel_{M(B)}$. 
Moreover, since $E^1(1) = 1$, $\parallel E^1\parallel =1$. Let $u_i\in A$ 
be an approximate unit of $B$. Then for any $b\in B, \mu\in\bf C$, 
$(b+\mu u_i)^*(b+\mu u_i) \geq 0$. Since $E$ is positive and of norm 1, 
$E(u_i^2) \leq 1$. Hence $E(b^*b)+\mid\mu\mid^2\cdot 1 + E(\overline\mu u_ib 
+\mu b^*u_i) \geq E((b+\mu u_i)^*(b+\mu u_i)) \geq 0$. Now the left hand 
side of the inequality converges to $E^1((b+\mu 1)^*(b+\mu 1))$ in norm. 
Hence, $E^1$ is positive. Now assume that $E^1((x+\mu 1)^*(x+\mu 1)) = 0$ 
and suppose that $\mu$ is non-zero. 
Then for any $a\in A$, $E((xa+\mu a)^*(xa+\mu a)) = 0$
and so $(xa+\mu a) = 0$ because $E$ is faithful. Let $p = E(-x/\mu)\in A$. 
Then $pa = a$ for all $a\in A$ which contradicts the assumption that $B$ 
is non-unital (see Lemma 1.1). Thus $\mu = 0$ and so $x=0$ as well (as $E$ 
is faithful). 

\par\medskip
Thus, we can construct ${\cal F}^1$ (the Hilbert $A^1$-module arising
from $E^1$; see the paragraph after Lemma 1.1) 
and the $C^*$-basic 
construction for $E^1$. Let $\lambda^1$ be the canonical map from 
$B^1$ to ${\cal L(F}^1)$ and $e^1_A$ be the element in ${\cal L(F}^1)$
that corresponds to $E^1$. 
Now since $B$ is a subalgebra of $B^1$ and $E^1$ 
extends $E$, $\cal F$ is a Hilbert $A^1$-submodule of ${\cal F}^1$ (note 
that $B\cdot A^1 \subseteq B$). 
It is clear that $e^1_A(\cal F) \subseteq F$ and 
its restriction to $\cal F$ is $e_A$. 
Let ${\cal K(F}^1)_0$ be the closed linear 
span in ${\cal L(F}^1)$ of the set $\{ \theta^1_{x,y}: x, y\in
\eta(B)\}$ (where $\theta^1_{x,y}(z) = x\langle y,z\rangle $ and 
let $\eta$ be the canonical map from $B^1$ to 
${\cal F}^1$). 
For any $T\in {\cal K(F}^1)_0$, $T(\cal F) \subseteq F$ (since it 
is true for each $\theta^1_{x,y}$). Define a map $\Psi$ from ${\cal K(F}^1)_0$ 
to $\cal L(F)$ by $\Psi(T) = T\mid_{\cal F}$ (the restriction of $T$). 
Then $\parallel \Psi \parallel \leq 1$ and $\Psi({\cal K(F}^1)_0) \subseteq 
\cal K(F)$. 
It is clear that $\Psi$ is a surjective *-homomorphism from 
${\cal K(F}^1)_0$ to $\cal K(F)$. 
Let $u_i$ be an approximate unit for $B$. Then for any $T\in {\cal K(F}^1)_0$, 
$T(\eta(u_i))$ will converge to $T(\eta(1))$ (since it is true for 
$\theta^1_{x,y}$ and $\parallel \eta(u_i)\parallel \leq 1$). Hence, if 
$\Psi(T)=0$, then $T(\eta(B))=0$ and so $T(\eta(B^1))=0$ which implies that 
$T=0$. Therefore, we obtain the following:

\par\medskip
\noindent {\bf Proposition 1.5:} Let ${\cal K(F}^1)_0$ be the closed linear 
span of the set $\{ \theta^1_{x,y}: x, y\in \eta(B)\}$ in ${\cal L(F}^1)$. 
Then ${\cal K(F}^1)_0 \cong \cal K(F)$. 
Consequently, $C^*\langle B,e_A\rangle $ is isomorphic 
to the closed linear span of the set $\{ \lambda^1(b) e^1_A \lambda^1(b'): 
b, b'\in B\}$ in $C^*\langle B^1, e^1_A\rangle $.

\par\medskip
We can now obtain the result that we want:

\par\medskip
\noindent {\bf Theorem 1.6:} Let B be a $C^*$-algebra and let 
$A$ be a $C^*$-subalgebra of $B$ that contains an approximate 
unit of $B$. Let $E$ be a faithful conditional expectation 
from $B$ to $A$. If $B$ acts on a Hilbert space $H$ 
faithfully (and non-degenerately) and $e$ is a projection on $H$ such that 

\noindent (1) $ebe = E(b)e$ for all $b\in B$ and 

\noindent (2) the map that sends $a\in A$ to $ae\in {\cal L}(H)$ 
is injective, 

\noindent then the norm closure of $BeB$ is canonically isomorphic to 
$C^*\langle B,e_A\rangle $.
\par
\noindent {\bf Proof:} If $B$ is unital, then it is just Theorem 1.3. 
Therefore, we assume that $B$ is non-unital. Since $B$ is represented 
faithfully on $H$, 
so is $M(B)$ and hence so is $B^1$. We would like to use Theorem 1.3. 
First of all, $e(b+\mu 1)e = E(b)e+\mu e = E^1(b+\mu 1)e$ and so condition (i)
of Theorem 1.3 holds. Secondly, since $A$ contains an approximate unit of 
$B$, $M(A)$ is represented non-degenerately on $H$ through a map $\phi$. 
Now let $K=eH$. Then, because of condition (1), 
condition (2) means that $A$ is represented faithfully and non-degenerately 
on $K$ (since $A$ contains an approximate unit of $B$) and hence 
so is $M(A)$. Now it is clear that this representation of $M(A)$ on $K$ is the 
restriction of $\phi$. Hence the map that sends $a\in A^1$ to $ae$ 
is injective. 
Now applying Theorem 1.3 to $(B^1, A^1, E^1)$, we know that the closure of 
$B^1eB^1$ is isomorphic to $C^*\langle B^1, e^1_A\rangle $ under an 
isomorphism $\psi$ which sends $beb'$ to $\lambda^1(b)e^1_A\lambda^1(b')$. 
Hence the restriction of $\psi$ 
will map the closure of $BeB$ to $\{ \lambda^1(b) e^1_A \lambda^1(b'): 
b, b'\in B\}$ which is isomorphic to $C^*\langle B,e_A\rangle $ by Proposition 1.5. 

\par\medskip\par\medskip
\noindent {\bf 2. Morita equivalences between fixed point algebras 
and crossed products}

\par\medskip
\par\medskip
In this section, we will prove a result about Morita equivalences
between fixed point algebras and crossed products for a special kind of
coaction of compact quantum groups (Theorem 2.13). 
We first prove the restricted case when $S=S_V$ (see [1]) for a regular 
multiplicative $V$ of compact type on a Hilbert space $H_V$ such that 
there exists a faithful Haar state $\varphi$ 
on $S$ and the coaction is injective and effective in the following sense. 

\par\medskip
\noindent {\bf Definition 2.2:}
A coaction $\epsilon$ on $B$ by $S$ is said to be effective if $\epsilon(B)
\cdot (B\otimes 1)$ is dense in $B\otimes S$. 
\par\medskip

Let $B$ be a $C^*$-algebra with an injective coaction $\epsilon$ by $S$ and 
let $B^{\epsilon}$ be the fixed point algebra. 
Let $E = (id\otimes \varphi)\epsilon$. 
It is not hard to show that $E$ is a faithful conditional expectation from 
$B$ to $B^{\epsilon}$. Moreover, $B^{\epsilon}$ contains an
approximate unit of $B$ by the following easy lemma:

\par\medskip
\noindent {\bf Lemma 2.1:} 
Let $B$ and $D$ be $C^*$-algebras and 
$\phi$ be a state on $D$. If $w_i$ is an approximate unit of $B\otimes D$, 
then $(id\otimes \phi)(w_i)$ is an approximate unit for $B$.

\par\medskip
Note that any dual coaction is effective since comultiplications are effective. 
From now on, we can use the same lines of proof as in [7, Section 5] to deduce 
Theorem 2.7. However, for completeness, we will repeat the arguments here. 

\par\medskip
\noindent {\bf Lemma 2.3:} $(id\otimes id\otimes \varphi)(V_{12}
V_{13}) = (id\otimes id\otimes \varphi)(V_{13})$. 
\par
\noindent {\bf Proof:} Since $\varphi$ is the Haar state, 
$(id \otimes \varphi)\delta((\omega\otimes id)V) = (\omega \otimes id
\otimes \varphi)(V_{13})$ for all $\omega\in {\cal L}(H_V)_*$. Now the 
lemma follows from the fact that ${\cal L}(H_V)_*$ separates points of 
${\cal L}(H_V)$. 
\par\medskip
\noindent {\bf Lemma 2.4:} Let $B$ be 
faithfully represented on a Hilbert space $H$. 
Regard $\epsilon$ as an injective map from $B$ to 
${\cal L}(H\otimes H_V)$ and let $e = 1\otimes p$ (where $p = 
(id\otimes \varphi)(V) \in {\cal L}(H_V))$. 
Then $\epsilon$ and $e$ will satisfy the 
two conditions of Theorem 1.6.
\par
\noindent {\bf Proof:} Since for any $a\in B^{\epsilon}$, 
$\epsilon (a) = a\otimes 1$, condition (2) of Theorem 1.6 is clear. 
For condition (1), observe that 
$(1\otimes p)\epsilon (b)(1\otimes p) = $ $(id\otimes id\otimes 
\varphi \otimes \varphi)(V_{23}(\epsilon (b)\otimes 1\otimes 1)
V_{24}) =$ $(id\otimes id\otimes \varphi \otimes \varphi)
(((id\otimes \delta_{V})\epsilon (b)\otimes 1)V_{23}V_{24})$ for all 
$b\in B$. 
Thus, using Lemma 2.3, $(1\otimes p)\epsilon (b)(1\otimes p) =$ 
$(id\otimes id\otimes \varphi \otimes \varphi)(((id\otimes \delta_{V})
\epsilon (b)\otimes 1)V_{24})=$ $[(id\otimes id\otimes \varphi)
((id\otimes \delta_{V})\epsilon (b))](1\otimes p)=$
$(id\otimes\varphi)(\epsilon (b))\otimes p$.
Hence we proved the lemma. 

\par\medskip
\noindent {\bf Lemma 2.5:} The set $P = 
\{ (id \otimes \varphi\cdot s)(V): s\in A_{V} \}$ is a dense 
subset of $\hat S_{V}$.
\par
\noindent {\bf Proof:} We first note that because $p\cdot \hat S_V = 
{\bf C}\cdot p$, $p\in \hat S_V$ (where $p=(id\otimes \varphi)(V)$). 
Moreover, if $s=(\omega\otimes id)(V)\in A_V$, 
then $(id\otimes \varphi)((1\otimes s)V)= 
(\omega\otimes id\otimes \varphi)(V_{13}V_{23}) = 
(id\otimes \varphi)(\omega\otimes id\otimes id)((\hat \delta_V\otimes id)V) = 
(\omega\otimes id)\hat \delta_V(p)$. Note that $\hat \delta_V(p)(x\otimes 1)
\in \hat S_V\otimes \hat S_V$ (for any $x\in \hat S_V$) 
and so $(\omega\otimes id)\hat \delta_V(p)\in 
\hat S_V$. Thus we have showen that $P$ is a subset of $\hat S_V$.
Let $t\in S_V$ be such that $(\varphi\cdot s)(t)=0$ for all $s\in 
A_V$. Then $\varphi(t^*t)=0$ (as $A_V$ is dense in $S_V$). 
Because $\varphi$ is faithful, $P' = \{ \varphi\cdot s: s\in A_V \}$ 
separates points of $S_V$. 
Hence $P'$ is $\sigma(S_V^*, S_V)$-dense 
in $S_V^*$. Therefore, for any $f\in S_V^*$, there exists a net 
${s_i}$ in $A_V$ such that $\varphi\cdot s_i$ converges to $f$ weakly. 
Note that $g(L_V(h)) = h(\rho_V(g))$ for all $h\in \hat S_V^*$ and 
$g\in S_V^*$ and that $L_V(\hat S_V^*)$ 
is a dense subset of $S_V$ (since $1\in S_V$). 
Hence for any $\nu \in {\cal L}(H_V)_*$, there exists a net $a_i$ in 
$P = \rho_V(P')$ such that $h(a_i)$ converges to $h(\rho_V(\nu))$
for any $h\in \hat S_V^*$. Therefore, the $\sigma (\hat S_V,\hat S_V^*)$-closure 
of $P$ will contain $\hat A_V$ and so $P$ is norm dense in 
$\hat S_V$ (because $P$ is a convex subset of $\hat S_V$). 

\par\medskip
\noindent {\bf Lemma 2.6:} 
Let the notation be the same as in Lemma 2.4. If, in addition, $\epsilon$ is 
effective, then the linear span, $T$, of 
$\{ \epsilon (a)(1\otimes p)\epsilon (b): a,b\in B\}$ is 
norm dense in $B\times_{\epsilon,r}\hat S_V$ 
\par
\noindent {\bf Proof:} 
We first note that $T$ is a subset of $B\times_{\epsilon,r}\hat S_V$. 
Since $\epsilon$ is a coaction, 
$(1\otimes p) \epsilon (b) = 
(id\otimes id\otimes \varphi)((\epsilon\otimes id)\epsilon(b)V_{23})$. 
Therefore, $\epsilon (a)(1\otimes p)\epsilon (b) = 
(id\otimes id\otimes \varphi)((\epsilon\otimes id)((a\otimes 1)\epsilon (b))
V_{23})$. Since $\epsilon$ is effective, elements of the form 
$(id\otimes id\otimes \varphi)
((\epsilon \otimes id)(c\otimes s)V_{23})$ ($c\in B$ and 
$s\in S_V$) can be approximated in norm by elements in $T$. Note that  
$(id\otimes id\otimes \varphi)((\epsilon \otimes id)(c\otimes s)V_{23})=$
$\epsilon (c)(1\otimes (id\otimes \varphi\cdot s)(V))$. Hence by Lemma 
2.5, $T$ is norm dense in $B\times_{\epsilon,r}\hat S_{V}$. 

\par\medskip
We can now state and prove the main theorem in this section.

\par\medskip
\noindent {\bf Theorem 2.7:} 
Let $B$ be a $C^*$-algebra with an injective 
and effective coaction $\epsilon$ by $S_V$. Let $B^{\epsilon}$ be the fixed 
point algebra of $\epsilon$. Then $B^{\epsilon}$ is strongly Morita 
equivalent to the reduced crossed product $B\times_{\epsilon,r}\hat S_V$
($=B\times_{\epsilon,max}\hat S_V$).
\par 
\noindent {\bf Proof:} By Lemma 2.4, Theorem 1.6 and the fact that the 
Hilbert-$B^{\epsilon}$-module $\cal F$ (defined by the conditional 
expectation) is full, $B^{\epsilon}$ 
is strongly Morita equivalent to the closure of the linear span of the set 
$\{ \epsilon (a)(1\otimes p)\epsilon (b): a,b\in B\}$ which, 
by Lemma 2.6, equals $B\times_{\epsilon,r}\hat S_V$.

\par\medskip
\par\medskip
We are now going to relax the conditions on $\epsilon$ and $S$. 
We first deal with the case when $\epsilon$ is not injective. 
As in [6, 2.17], we consider $I=Ker(\epsilon)$. Then there is a coaction 
$\epsilon '$ on $A=B/I$ given by $\epsilon '(q(b))=(q\otimes id)\epsilon(b)$ 
for any $b\in B$ (where $q$ is the canonical quotient from $B$ to $A$). 
It is easily seen that if $\epsilon$ is effective, then so is $\epsilon '$. 
Note that $(id\otimes \varphi)\epsilon(B) = B^{\epsilon}$ no matter $\epsilon$ 
is injective or not.
We first show the following lemma:

\par\medskip
\noindent {\bf Lemma 2.8:} 
With the notation above, $A^{\epsilon '}$ is isomorphic to $B^{\epsilon}$. 
\par
\noindent {\bf Proof:} Note that $A^{\epsilon '} = (id\otimes \varphi)
\epsilon '(q(B)) = (q\otimes \varphi)\epsilon(B) = q(B^{\epsilon})$.
Now if $b\in B^{\epsilon}$ be such that $q(b)=0$, then $b\in I=Ker(\epsilon)$ 
and $\epsilon (b)=b\otimes 1$ which implies that $b=0$. Hence, the restriction 
of $q$ on $B^{\epsilon}$ is injective and is therefore an isomorphism. 

\par\medskip
\noindent {\bf Proposition 2.9:} 
Let everything be the same as Theorem 2.7 except that we don't assume 
$\epsilon$ to be injective. Then $B^{\epsilon}$ is strongly Morita 
equivalent to the reduced crossed product $B\times_{\epsilon,r}\hat S_V$.
\par
\noindent {\bf Proof:} 
Since $V$ is of compact type, it is automatically amenable. Therefore, by 
[6, 2.19(a)], we have $A\times_{\epsilon ',r}\hat S_V = B\times_{\epsilon,r}
\hat S_V$ (where $A$ and $\epsilon '$ be as defined above). 
Therefore, $B^{\epsilon}$ is strongly Morita equivalent to 
$B\times_{\epsilon,r}\hat S_V$ by Theorem 2.7.

\par\medskip
We now turn to the case when $S=(S_V)_p$ (see [6]). 

\par\medskip
\noindent {\bf Lemma 2.10:} Let $\varphi_p = \varphi\circ L_V$.
Then $\varphi_p$ is an invariant state on $(S_V)_p$ (in the sense that 
$(id\otimes \varphi_p)\delta_p = \varphi_p\cdot 1$).
\par
\noindent {\bf Proof:} 
Let $p$ be the minimum central projection in $\hat S_V$ that 
corresponds to $\varphi$. 
Let $\chi$ be the injective norm decreasing algebraic homomorphism from 
$(S_V)_p^*$ to $M(\hat S_V)$ as given in [6, A6]. Then it is clear that
$\chi(\varphi_p) = p$ and hence $\varphi_p$ is a central minimum projection 
in $(S_V)_p^*$. This means that $(f\otimes \varphi_p)\delta_p = 
f(1)\varphi_p$ for any $f\in (S_V)_p^*$ (note that $\delta_p(1)=1\otimes 1$). 

\par\medskip
Let $B$ be a $C^*$-algebra with coaction $\epsilon ''$ by $(S_V)_p$.
As in [6, 2.14], we consider $\epsilon = (id\otimes L_V)\epsilon ''$. 
It is clear that if $\epsilon ''$ is effective, then so is $\epsilon$. 
Note that we also have $(id\otimes \varphi_p)\epsilon ''(B) = B^{\epsilon ''}$. 
Now $B^{\epsilon ''} = (id\otimes \varphi_p)\epsilon ''(B) = 
(id\otimes \varphi)\epsilon(B) = B^{\epsilon}$.
Moreover, we recall from [6, 2.14] that $B\times_{\epsilon '',r}(\hat S_V)_p$ 
is isomorphic to $B\times_{\epsilon,r}\hat S_V$. These, together with 
Proposition 2.9, proved the following generalisation of Theorem 2.7:

\par\medskip
\noindent {\bf Proposition 2.11:}
Let $B$ be a $C^*$-algebra with an 
effective coaction $\epsilon$ by $S = S_V$ or $(S_V)_p$. 
Let $B^{\epsilon}$ be the fixed point algebra of $\epsilon$. 
Then $B^{\epsilon}$ is strongly Morita equivalent to the reduced 
crossed product $B\times_{\epsilon,r}\hat S$ 
($=B\times_{\epsilon,max}\hat S$). 

\par\medskip
Now we would like to further generalise $S$. Let $S$ be a compact quantum 
group (i.e. an unital Hopf $C^*$-algebra with comultiplication $\delta$ 
such that both $\delta (S)(S\otimes 1)$ and $\delta (S)(1\otimes S)$ 
are dense in $S\otimes S$). 
By [11], $S$ has a Haar state $\phi$. Let $(H, \pi)$ be 
the GNS representation corresponding to $\phi$ and $V$ be the multiplicative 
unitary on $H$ as defined in [1, 1.2(4)]. Then $\pi$ is a surjective map from 
$S$ to $S_V$. Let $\varphi$ be the Haar state on $S_V$ (which equals the 
state defined by $\phi$). Then $\varphi$ is faithful by the following lemma.

\par\medskip
\noindent {\bf Lemma 2.12:} 
Let $A$ be a compact quantum group with a Haar state $\psi$. Then $\psi$ 
is a faithful state if and only if the GNS representation $(H, \pi)$
corresponding to $\psi$ is faithful.

\noindent {\bf Proof:} 
Suppose that $(H, \pi)$ is faithful. Let $N=\{ x\in A: \psi (x^*x)=0 \}$. 
Then it is clear that $N=\{ x\in A: \psi (yx)=0$ for all $y\in A \}$. By 
[15, 5.6(6)], $\psi (ab) = \psi(b (f_1*a*f_1))$ for all $a\in {\cal A}$. 
Now since ${\cal A}$ is a dense subalgebra of $A$ and $f_1*{\cal A}*f_1 = 
{\cal A}$, $N=\{ x\in A: \psi (xy)=0$ for all $y\in A \}$. Hence $N$ is an 
ideal of $A$ and so $N = Ker(\pi) = 0$. Thus $\psi$ is faithful. 

\par\medskip
Let $B$ be 
a $C^*$-algebra with coaction $\epsilon ''$ by $S$ and let 
$\epsilon = (id\otimes \pi)\epsilon ''$. Then $\epsilon$ is a coaction 
on $B$ by $S_V$. By the same argument as in Proposition 2.11, 
we have $B^{\epsilon} = 
B^{\epsilon ''}$. Therefore, if we define the reduced crossed product 
of $\epsilon ''$ to be the reduced crossed product of $\epsilon$, then 
we have the following:

\par\medskip
\noindent {\bf Theorem 2.13:}
Let $B$ be a $C^*$-algebra with an 
effective coaction $\epsilon$ by a compact quantum group $S$. 
Then $B^{\epsilon}$ is strongly Morita equivalent to the reduced 
crossed product $B\times_{\epsilon,r}\hat S$. 

\par\medskip
Note that the assumption of effectiveness cannot be removed. For 
example, if we consider $\epsilon$ to be the trivial coaction on $B$ (i.e. 
$\epsilon(b) = b\otimes 1$), then the fixed point algebra is $B$ itself 
while the reduced crossed product is $B\otimes \hat S$ which are clearly 
not Morita equivalent (e.g. when $B=\bf C$ and $S$ is non-trivial). 

\par\medskip
\noindent {\bf Remark 2.14:} 
(a) We was told that van Daele and Zhang have recently proved a similar 
result to Theorem 2.13 (but in a purely algebraic framework) for 
algebraic quantum groups (which can be regarded as a generalization of
both compact quantum groups and discrete quantum groups). 
For a reference, we refer the readers to [12]. 
\par\noindent
(b) By the proof of Theorem 2.7, we know that in general (i.e. when $\epsilon$ is 
not effective), $B^{\epsilon}$ is strongly Morita equivalent to the ideal 
$I = \{ \epsilon(a) (1\otimes p) \epsilon(b): a,b \in B\}$ of 
$B\times_{\epsilon,r}\hat S_V$ ($I$ is an ideal since $p$ is a central 
minimum projection). Consequently, if the reduced crossed product satisfies 
some properties that are preserved under taking ideal and under Morita 
equivalence, then so does the fixed point algebra. In particular, we have 
the following corollary.

\par\medskip
\noindent {\bf Corollary 2.15:}
Let $B$ be a $C^*$-algebra with a coaction (not necessary effective) 
$\epsilon$ by a compact quantum group $S$. Then 
\par
\noindent (a) if $B\times_{\epsilon,r}\hat S$ is nuclear (this is the case 
if $B$ is, see [6, 3.4]), so is $B^{\epsilon}$;
\par
\noindent (b) if $B\times_{\epsilon,r}\hat S$ is liminal (resp. postliminal), 
then so is $B^{\epsilon}$; 
\par
\noindent (c) if $B\times_{\epsilon,r}\hat S$ is simple, then it is 
strongly Morita equivalent to $B^{\epsilon}$ and so 
$B^{\epsilon}$ is simple.

\par\medskip
Finally, we would like to apply Theorem 2.7 to the case of coactions by 
discrete groups. We recall from [4, 2.6], that if $A$ is a $C^*$-algebra 
with ``reduced'' coaction $\epsilon$ by a discrete group $G$ (i.e. 
an injective and non-degenerate coaction by $C^*_r(G)$), then $A = 
\overline{\bigoplus_{t\in G} A_t}$ (with $A_e$ being the fixed point 
algebra of $\epsilon$). 
Let $\psi$ be the Haar state on $C^*_r(G)$ and $\lambda_t$ be the 
canonical image of $t\in G$ in $C_r^*(G)$.
Let $\psi_t = \psi\cdot \lambda_{t^{-1}}$. 
Then $A_t = (id\otimes \psi_t)\epsilon(A)$. 
Recall from [5, Section 5] that $\epsilon$ is said to 
be full if $\overline{A_t\cdot A_{t^{-1}}} = A_e$ for all $t\in G$ (or 
equivalently, $\overline{A_r\cdot A_s} = A_{rs}$ for all $r,s \in G$). 

\par\medskip
\noindent {\bf Lemma 2.16:} $\epsilon$ is effective if and only if 
it is full. 
\par
\noindent {\bf Proof:} 
Suppose that $\epsilon$ is effective. Then by acting $(id\otimes \psi_t)$ on 
both sides of the equation $\overline{\epsilon(A)\cdot (A\otimes 1)} = 
A\otimes C_r^*(G)$, we have $\overline{A_t\cdot A} = A$. 
Thus by applying $(id\otimes \psi_{-t})\epsilon$ again, we obtain 
$\overline{A_t\cdot A_{t^{-1}}} = A_e$. 
Conversely, suppose that $\epsilon$ is full.
It is required to show that 
$\overline{\epsilon(A)\cdot (A\otimes 1)} \supseteq 
A_r\otimes \lambda_s$. Actually, it suffices to show that 
$(A_s\otimes \lambda_s)\cdot (A_{s^{-1}r}\otimes 1) \supseteq 
A_r\otimes \lambda_s$. However, it is clear from the definition of full 
coaction that this relation holds.

\par\medskip
Now using Theorem 2.7, we can obtain as a corollary [5, 5.2] 
(which stated that the fixed point algebra and the crossed product 
are Morita equivalent if the coaction is full). 

\par\medskip\par\medskip
\noindent {\bf 3. An application of the main theorem}

\par\medskip
\par\medskip

In this section, we will use the main result in Section 2 to prove a
imprimitivity type theorem for crossed products of coactions by
discrete quantum groups. 
We first recall some basic definition from [7, Section 4].

\par\medskip
\noindent {\bf Definition 3.1:}
Let $U$ and $V$ be regular multiplicative unitaries on the Hilbert spaces 
$H$ and $K$ respectively. $X\in {\cal L}(K\otimes H)$ 
is said to be a $U$-$V$-birepresentation if $X$ is a representation
of $V$ as well as a corepresentation of $U$.

\par\medskip
We recall from [7, 3.9] that a $U$-$V$-birepresentation $X$ will induce 
a Hopf *-homomorphism $L_{X''}$ from $(S_{U})_{p}$ to 
$M[(S_{V})_{p}]$ as well as a Hopf *-homomorphism $\rho_{X'}$
from $(\hat{S}_{V})_{p}$ to $M[(\hat{S}_{U})_{p}]$. Therefore, we have a 
coaction $\epsilon_X = (id\otimes \rho_{X'})\hat\delta_V$ 
on $(\hat{S}_{V})_{p}$ by $(\hat{S}_{U})_{p}$. 

\par\medskip
\noindent {\bf Definition 3.2:}
Let $U$, $V$, $W$ be regular multiplicative unitaries. 
\par
\noindent (a) $W$ is said to be a sub-multiplicative unitary of $V$ 
if there exists a $V$-$W$-birepresentation $Y$ such that 
$L_{Y^{\prime\prime}}((S_{V})_{p}) = (S_{W})_{p}$.
\par
\noindent (b) $U$ is said to be a quotient of $V$ if there exists a 
$U$-$V$-birepresentation
$X$ such that $\rho _{X^\prime }((\hat{S}_{V})_{p}) = (\hat{S}_{U})_{p}$.
\par
\noindent (c) Let $U$, $V$ and $W$ be of discrete type such 
that $U$ is a quotient of $V$ through $X$ and 
$W$ is a submultiplicative unitary of $V$ through $Y$. 
Then $W$ is said to be normal if $\rho_{Y'}$ is an isomorphism from 
$(\hat{S}_{W})_{p}$ to the fixed
point algebra of $\hat{\epsilon}_{X}$ in $(\hat{S}_{V})_{p}$.

\par\medskip
The above terminologies come from the case of locally compact groups. 
If $H$ and $G$ are two locally compact groups, then $H$ is a quotient 
of $G$ if and only if $C^*(H)$ is a quotient of $C^*(G)$ as Hopf 
$C^*$-algebras. 
Moreover, $H$ is a subgroup of $G$ if and only if $C_0(H)$ is a quotient 
of $C_0(G)$ (as Hopf $C^*$-algebras).
A normal subgroup of a discrete group will certainly satisfy condition 
3.2(c) (see [7,A2]). 
Using Proposition 2.11, we can now give a positive answer to the 
remark in [7, 5.7(a)] (i.e. we don't need to assume the amenability 
for $U$ in the imprimitivity theorem). More precisely, we have:

\par\medskip
\noindent {\bf Theorem 3.3:} 
Let $U$, $V$ and $W$ be regular multiplicative unitaries of discrete type 
such that $U$ comes from a discrete Kac algebra. If $W$ is a normal 
submultiplicative unitary of $V$ with quotient $U$, then $(\hat S_W)_p$ 
is strongly Morita equivalent to $(\hat S_V)_p\times_{\epsilon',r}S_U$.

\par\medskip
Moreover, we can have the following version of the imprimitivity theorem for 
crossed products: Let $U$, $V$ and $W$ be multiplicative unitaries that 
come from discrete Kac algebras such that $W$ is a normal sub-multiplicative 
unitary of $V$ (through a $V$-$W$-bi-representation $Y$) with quotient $U$
(through a $U$-$V$-bi-representation $X$). 
Furthermore, suppose that $W$ is amenable 
(note that all $U$, $V$ and $W$ are co-amenable in this case). 

\par\medskip
\noindent {\bf Theorem 3.4:} 
Let $U$, $V$ and $W$ be as given in the previous paragraph.  
Let $A$ be a $C^*$-algebra with an injective and non-degenerate 
coaction $\epsilon$ by $S_V = (S_V)_p$ and let $\hat \epsilon$ be 
the dual coaction on the full crossed product 
$A\times_{\epsilon}\hat S_V$ (see [7, 1.12]). 
Suppose that $\epsilon '$ and $\epsilon ''$ are the coactions on 
$A$ and $A\times_{\epsilon}\hat S_V$ by 
$S_W$ and $(\hat S_U)_p$ induced from $\epsilon$ and $\hat \epsilon$ 
respectively (i.e. $\epsilon' = (id\otimes L_{Y''}\epsilon$ and 
$\epsilon '' = (id\otimes \rho_{X'})\hat\epsilon$). 
Then $A\times_{\epsilon '}\hat S_W$ is strongly 
Morita equivalent to $(A\times_{\epsilon}\hat S_V)\times_{\epsilon '',r}S_U$. 

\par\medskip
We need several lemmas to prove this theorem. 

\par\medskip
\noindent {\bf Lemma 3.5:} 
If $\epsilon$ is an injective and non-degenerate coaction on a 
$C^*$-algebra $A$ by $S_V$, then the induced coaction $\epsilon ' = 
(id\otimes L_{Y''})\epsilon$ on $A$ by $S_W$ is also injective 
and non-degenerate. 
\par
\noindent {\bf Proof:} 
Let $E_V$ be the co-identity on $S_V = (S_V)_p$. Then for any $a\in A$, 
$\epsilon((id\otimes E_V)\epsilon(a)-a) = 0$ and so 
$(id\otimes E_V)\epsilon(a)=a$. 
Now if $E_W$ is the co-identity on $S_W$, then $E_V=E_W\circ L_Y$ by 
[7, 3.2(d)]. Therefore, $(id\otimes E_W)\epsilon '(a) = a$ and 
so $\epsilon '$ is injective. 
$\epsilon '$ is non-degenerate since $L_Y$ is surjective. 

\par\medskip
\noindent {\bf Lemma 3.6:}
Let $B$ and $C$ be $C^*$-algebras and $\phi$ a non-degenerate 
*-homomorphism from $B$ to $C$ with kernel $I$. Let $\overline\phi$ 
be the extension of $\phi$ from $M(B)$ to $M(C)$. 
Then $Ker(\overline\phi) = \{ m\in M(B): m\cdot B, B\cdot m \subseteq I\}$. 
\par
\noindent {\bf Proof:} 
Let $Q$ be the quotient map from $B$ to $B/I$ and $\psi$ be the canonical 
non-degenerate monomorphism from $B/I$ to $C$. Then $\overline\phi = 
\overline\psi \circ \overline{Q}$. Now the lemma follows from the easy 
facts that $Ker(\overline{Q})=\{ m\in M(B): m\cdot B, B\cdot m \subseteq I\}$
and that $\overline\psi$ is injective. 

\par\medskip
\noindent {\bf Lemma 3.7:} 
Let $B$ and $C$ be two $C^*$-algebras with coactions $\epsilon_B$ and 
$\epsilon_C$ respectively by a Hopf $C^*$-algebra $S$. Suppose that $S$ 
is $C^*$-exact (as a $C^*$-algebra). If $\psi$ is a non-degenerate 
equivariant *-homomorphism 
from $B$ to $C$, then $Ker(\psi)$ is a weakly invariant ideal of $A$ (in the 
sense of [8, 3.14]). 
\par
\noindent {\bf Proof:} 
Let $I=Ker(\psi)$ and $q$ be the canonical quotient from $B$ to $B/I$. 
Since $S$ is $C^*$-exact, $Ker(q\otimes id) = I\otimes S$. 
Hence $Ker(\psi\otimes id) = Ker (\tilde\psi\circ q\otimes id) = 
I\otimes S$ (where $\tilde\psi$ is the canonical injection from $B/I$ to $C$). 
Therefore, by Lemma 3.6, 
$Ker(\overline{\psi\otimes id}) = 
\{ m\in M(B\otimes S): m\cdot (B\otimes S), (B\otimes S)\cdot m 
\subseteq I\otimes S\}$. 
Now for any $x\in I$, $(\overline{\psi\otimes id})
\epsilon_B(x) = \epsilon_C(\psi(x)) = 0$. Thus, $\epsilon_B(x)\in 
\tilde M(B\otimes S)\cap Ker(\overline{\psi\otimes id})\subseteq 
\tilde M(I\otimes S)$ (by [2, 1.4]). 

\par\medskip

Let $D$ be a $C^*$-algebra with an injective coaction $\delta$ by a Hopf 
$C^*$-algebra that is defined by a regular multiplicative unitary. Let 
$(C, j_D, \nu)$ be the full crossed product of $\delta$ (see [6]). 
Then $j_D$ is injective because the canonical map from $D$ to the 
reduced crossed product is injective and the reduced representation 
is covariant. 
\par
Stimulated by [9, 4.3], we have the following lemma.

\par\medskip
\noindent {\bf Lemma 3.8:} 
Let $B$ be a $C^*$-algebra with an injective and non-degenerate 
coaction $\epsilon$ 
by $S=S_{V'}$ where $V'$ is a co-amenable irreducible multiplicative 
unitary on a Hilbert space $H$. 
Let $D$ be a $C^*$-algebra with an injective 
coaction $\delta$ by $\hat S = \hat S_{V'}$. 
If $\phi$ is an (non-degenerate) equivariant *-homomorphism from 
$B\times_{\epsilon,r}\hat S$ to $D$ such that the restriction of $\phi$ 
on $B$ (considered as a subalgebra of $M(B\times_{\epsilon,r}\hat S)$) 
is injective, then $\phi$ is injective.
\par
\noindent {\bf Proof:} 
Since $\phi$ is equivariant, by [6, 3.9], it induces a map $\Phi$ from 
$(B\times_{\epsilon,r}\hat S)\times_{\hat\epsilon} S$ to $D\times_{\delta} S$ 
such that $\Phi\circ j = j_D\circ \phi$ and $\Phi\circ \mu =\nu$ 
(where $((B\times_{\epsilon,r}\hat S)\times_{\hat\epsilon} S, j, \mu)$ and 
$(D\times_{\delta} S, j_D, \nu)$ are the full crossed products of 
$\hat\epsilon$ and $\delta$ respectively). 
It is not hard to show that $\Phi$ is equivariant with respect to the 
dual coactions defined in [7, 1.12]. 
Now since $V'$ is a co-amenable irreducible multiplicative unitary, 
$(B\times_{\epsilon,r}\hat S)\times_{\hat\epsilon} S = 
(B\times_{\epsilon,r}\hat S)\times_{\hat\epsilon,r} S \cong 
B\otimes {\cal K}(H)$ (by [1, 7.5]). 
It is also clear that the dual coactions defined in [7, 1.12] and 
in [1, Section 7] are identical in this case. 
Now by the proof of [1, 7.5], for any $b\in B$, $\pi_L(b)\otimes 1$ 
($\in B\times_{\epsilon,r}\hat S\times_{\hat\epsilon,r} S$) 
corresponds to $\pi_R(b) = (id\otimes R)\epsilon(b)$ 
($\in B\otimes {\cal K}(H)$) under the isomorphism. 
Regard $\Phi$ as a map from $B\otimes {\cal K}(H)$ to 
$D\times_{\delta} S$. 
Then we have $\Phi(\pi_R(b)) = j_D(\phi(b))$. 
Let $J = Ker(\Phi)$. 
$J$ is weakly invariant by Lemma 3.7 (note 
that $S$ is nuclear by [6, 3.6]). Moreover, by 
the proof of [8, 3.15], $J = I\otimes {\cal K}(H)$ for some weakly invariant 
ideal $I$ of $B$. Now for any $x\in I$, $\epsilon(x)\in 
\tilde M(I\otimes S)$. Since $R$ is an isomorphism from 
$S$ to $\hat S_{\tilde V} \subseteq {\cal L}(H)$,  $\pi_R(x)\in 
\tilde M(I\otimes \hat S_{\tilde V}) \subseteq Ker(\overline\Phi)$ 
(by Lemma 3.6). Hence $j_D(\phi(x)) = \Phi(\pi_R(x)) = 0$ and so $\phi(x)=0$
(note that $j_D$ is injective as $\delta$ is injective). 
Therefore, from the hypothesis, $I=0$ and so $\Phi$ is injective. Using 
the facts that $j$ is injective (the dual coaction $\hat \epsilon$ is 
injective) and that $\Phi\circ j = j_D\circ \phi$, we showed that $\phi$ 
is also injective.

\par\medskip
Now we can give a proof for Theorem 3.4. 

\par\medskip
\noindent {\bf Proof:} (Theorem 3.4) 
Let $(A\times_{\epsilon}\hat S_V, j, \mu)$ be the full crossed product of 
$\epsilon$. Then $(j\otimes id)\epsilon '(a) = 
(id\otimes L_Y) ((\mu\otimes id)(V')(j(a)\otimes 1)(\mu\otimes id)(V')^*) = 
((\mu\otimes id)(Y')(j(a)\otimes 1)(\mu\otimes id)(Y')^*)$. Hence, 
$(j, \mu\circ \rho_{Y'})$ is a covariant pair for $\epsilon '$ and so there 
exists a non-degenerate *-homomorphism $\phi$ from $A\times_{\epsilon '}\hat S_W$ 
to $M(A\times_{\epsilon}\hat S_V)$. Moreover, since $1\in (\hat S_V)_p$, 
$j(A) \subseteq A\times_{\epsilon}\hat S_V$. Hence $\phi 
(A\times_{\epsilon '}\hat S_W) \subseteq A\times_{\epsilon}\hat S_V$. 
Let $\varphi$ be the Haar state on $(\hat S_U)_p$ and 
let $D = (A\times_{\epsilon}\hat S_V)^
{\epsilon ''}$. Then $D = (id \otimes \varphi)\epsilon ''
(A\times_{\epsilon}\hat S_V)$. 
Now for any $a\in A$ and $t\in (\hat S_V)_p$, $(id \otimes \varphi)\epsilon ''
(j(a)\mu(t)) = (id \otimes \varphi\circ \rho_{X'})((j(a)\otimes 1)
(\mu \otimes id)\hat\delta_V(t)) = j(a)(\mu \otimes \varphi)\epsilon_{X'}(t)$. 
By the assumption that $\rho_{Y'}$ is an isomorphism from 
$(\hat S_W)_p$ to $(\hat S_V)_p^{\epsilon_{X'}}$, we have 
$\phi (A\times_{\epsilon '}\hat S_W) = D$. 
Next, we would like to show that $\phi$ is injective. 
Note that for any $a\in A$ and $s\in (\hat S_W)_p$, $\hat\epsilon(j(a)\mu
(\rho_{Y'}(s))) = (j(a)\otimes 1)(\mu\otimes id)(\rho_{Y'}\otimes \rho_{Y'})
(\hat\delta_W(s)) \in (A\times_{\epsilon}\hat S_V)^{\epsilon ''}\otimes 
(\hat S_W)_p$ (since $(\hat S_W)_p$ is unital and $\rho_{Y'}$ is a Hopf 
*-monomorphism). Moreover, if $u_i$ is an 
approximate unit of $A$, then $\hat\epsilon(j(u_i)) = j(u_i)\otimes 1$ is 
an approximate unit for $(A\times_{\epsilon}\hat S_V)^{\epsilon ''}\otimes 
(\hat S_W)_p$. Hence $\hat\epsilon$ induces an injective coaction $\delta$ on 
$(A\times_{\epsilon}\hat S_V)^{\epsilon ''}$ by $(\hat S_W)_p = \hat S_W$ 
(note that $\delta$ is injective since $(id\otimes \hat E_W)\delta = id$ 
with $\hat E_W$ being the co-identity of $(\hat S_W)_p$). 
Furthermore, $\phi$ is clearly 
injective on $A$ since $j$ is injective (note that $\epsilon$ is injective). 
Therefore, by Lemma 3.8, $\phi$ is injective. 
We are now going to use Proposition 2.11 to prove this theorem. 
It remains to show 
that $\epsilon ''$ is effective. Since $\rho_{X'}$ is surjective 
(by assumption), it suffices to show that $\hat\epsilon$ is effective. 
However, since $\hat\delta_V((\hat S_V)_p)\cdot ((\hat S_V)_p\otimes 1) = 
(\hat S_V)_p\otimes (\hat S_V)_p$, it is easily seen that $\hat\epsilon$ is 
effective. 
\par\medskip
\par\medskip

\noindent {\it Acknowledgement}
\par\medskip
This work was done during our visit to the Fields Institute. We would like 
to thank the organiser of the special program for operator algebras 94-95, 
Prof. G. Elliott, for his hospitality during our visit. 

\par
\medskip
\par
\medskip
\par
\medskip
\noindent {\bf REFERENCE:}
\par
\medskip
\noindent [1] S. Baaj and G. Skandalis, Unitaires multiplicatifs et 
dualit\'e pour les
produits crois\'es de $C^{*}$-alg\`ebres, Ann. scient. \'Ec. Norm. Sup., 
$4^{e}$ s\'erie, t. 26
(1993), 425-488.
\par
\noindent [2] M. B. Landstad, J. Phillips, I. Raeburn and C. E. Sutherland, 
Representations of crossed products by coactions and principal bundles, 
Tran. of Amer. Math. Soc., vol 299, no. 2 (1987), 747-784.
\par
\noindent [3] K. Mansfield, Induced representations of crossed products by 
coactions, J.
Funct. Anal., 97 (1991), no. 1, 112-161.
\par
\noindent [4] C. K. Ng, Discrete coactions on $C^*$-algebras, Journal of 
Austral. Math. Soc. (Series A), 60 (1996), 118-127.
\par
\noindent [5] C. K. Ng, Discrete coactions on Hilbert $C^*$-modules, Math. 
Proc. of Camb. Phil. Soc., 119 (1996), 103-112.
\par
\noindent [6] C. K. Ng, Coactions and crossed products of Hopf 
$C^{*}$-algebras, Proc. London Math. Soc.(3), 72 (1996), 638-656.
\par
\noindent [7] C. K. Ng, Morphisms of regular multiplicative unitaries, 
Journal of Operator Theory, to appear.
\par
\noindent [8] C. K. Ng, Coactions on Hilbert $C^*$-modules, preprint.
\par
\noindent [9] I. Raeburn, On crossed products by coactions and their 
representation theory, Proc. London Math. Soc. (3), 64 (1992), 625-652.
\par
\noindent [10] M. A. Rieffel, Unitary representations of group extensions: 
An algebraic approach to the theory of Mackey and Blattner, in 
``Adv. in Math. Suppl. Stud.'', Vol 4, pp. 43-82, Academic Press, Orlando, 
FL, 1979.
\par
\noindent [11] A. van Daele, The Haar Measure on a compact quantum group, 
preprint.
\par
\noindent [12] A. van Daele and Y. Zhang, Galois theory for multiplier 
Hopf algebras with integrals, preprint.
\par
\noindent [13] S. Z. Wang, Free products of compact quantum groups, 
Comm. Math. Phys., 167 (1995), 671-692.
\par
\noindent [14] Y. Watatani, Index for $C^*$-algebras, Memoirs of Amer. Math. 
Soc. no. 242 (1990).
\par
\noindent [15] S. L. Woronowicz, Compact matrix pseudogroups, Comm. Math. 
Phys. 111(1987), 613-665.
\par
\medskip
\noindent MATHEMATICAL INSTITUTE, OXFORD UNIVERSITY, 24-29 ST. GILES, OXFORD 
OX1 3LB, UNITED KINGDOM.
\par
\noindent $E$-mail address: ng@maths.ox.ac.uk
\par
\end{document}